\documentclass[12pt,thmsa]{article}
\usepackage{amsfonts}
\usepackage{amssymb}
\usepackage{graphicx}
\usepackage{amssymb,amsmath}
\newtheorem{teo}{Theorem}[section]

\newtheorem{lema}[teo]{Lemma}

\newtheorem{coro}[teo]{Corollary}
\newtheorem{lemma}[teo]{Lemma}

\newtheorem{obs2}[teo]{Remark}

\newtheorem{tea}{Theorem}[subsection]

\newtheorem{no2}[teo]{Note}

\newtheorem{no3}[tea]{Note}

\newcommand{\Gal}{{\rm Gal}}
\newcommand{\Frob}{{\rm Frob }}

\newcommand{\PSL}{{\rm PSL}}
\newcommand{\SL}{{\rm SL}}

\newcommand{\PGL}{{\rm PGL}}

\newcommand{\GL}{{\rm GL}}

\newcommand{\F}{{\mathbb{F}}}
\newcommand{\Q}{{\mathbb{Q}}}

\def\timehm{\count31=\time \count32=\count31 \divide\count31 by 60
\number\count31 \multiply\count31 by 60 \advance\count32 by
-\count31 :\ifnum\count32<10 0\fi \number\count32}

\newcommand{\qed}{\hfill\rule{2mm}{2mm}}

\def\mod{\mathop{\rm mod}\nolimits}

\def\ideal#1{<\kern-2pt #1\kern-2pt >}

\begin{document}
\title{{\bf  Langlands Base Change for $\GL(2)$}
}
\author{Luis Dieulefait\thanks{Research partially supported
by MICINN grant MTM2009-07024 and by an ICREA Academia Research Prize}
\\
Universitat de Barcelona\\
e-mail: ldieulefait@ub.edu\\
 }
\date{\empty}

\maketitle

\vskip -20mm

\begin{abstract} Let $F$ be a totally real Galois number field. We prove the existence of base change relative to the extension $F/\Q$ for every holomorphic newform of weight at least $2$ and odd level, under simple local assumptions on the field $F$. 
\end{abstract}

\section{Introduction}
In their 1997 paper [HM], Hida and Maeda proposed a strategy to attack the problem of non-abelian base change for a totally real extension $F$. The case of solvable base change was known to be true by the work of Langlands (cf. [L]). Given a newform $f:= f_1$, they propose to find a sequence of links (congruences modulo suitable primes) starting from $f$ and ending in some newform $f_j$ (and calling $f_2$ to $f_{j-1}$ the modular forms appearing as \it{cha\^{i}nons} \rm of this chain) such that for some reason (in their case, they take a CM form) it is known that $f_j$ can be lifted to a Hilbert modular form on $F$. Then, assuming that for the restrictions to the absolute Galois group $G_F$ of $F$ of the Galois representations in this chain suitable Modularity Lifting Theorems (M.L.T.) apply at all the links, we can propagate modularity (over $F$) from the restriction to $G_F$ of $f_{i+1}$ to that of $f_{i}$ ($i= j-1, j-2,....,1$), thus proving that $f$ can be lifted to $F$. In other words, through suitable congruences and M.L.T., the liftability of $f_j$ to $F$ implies that of $f_{j-1}$, and so on until deducing the liftability of the given $f$, the first newform in the chain of congruences. They manage this way to prove that base change holds for infinitely many modular forms of prime power level.\\
In the new decade new M.L.T. (over totally real number fields) have been proved, and in particular recent results of Kisin are strong enough to suggest that this strategy can now be applied to prove base change in almost full generality. Recent M.L.T. still have conditions on the size of the residual image and on its restriction to the decomposition group at $p$, but on the other hand in the recent proofs of Serre's conjecture (cf. [S]) given by the author and by Khare-Wintenberger (see [D] and [KW1], [KW2]) new \it{astuces} \rm have been developed to propagate modularity despite these conditions. In fact, combining the Hida-Maeda strategy with strong Modularity Lifting Theorems of Kisin (cf. [KiFM], [KiBT] and [Ki2]), Skinner and Wiles (cf. [SW]) and Geraghty (cf. [G]), the \it{propagatory techniques} \rm (of the author and K.-W.) just mentioned, and some new ideas, we can prove base change for $\GL(2)$ under some mild assumptions, namely, the following is true:

\begin{teo} \label{teo:main}
Let $F$ be a totally real Galois number field. Let $f$ be a holomorphic newform of weight at least $2$ and odd level $N$. Assume that the following two conditions are satisfied:
\begin{enumerate}
	\item the primes $2,3, 7$ and $11$ are split in $F$.
	\item if $5 \mid N$, then $5$ is split in $F$. 
\end{enumerate}
Then,  $f$ is liftable to $F$, i.e., there is a Hilbert modular form $\hat{f}$ over $F$ such that the restrictions of the $\lambda$-adic Galois representations attached to $f$ to $G_F$ agree with the Galois representations attached to $\hat{f}$.
\end{teo}

We will also include in the last section some elementary corollaries of our base change result.\\
\newline
Remark: for modular forms that can also be obtained from a definite quaternion algebra, assuming that the Galois group of $F/\Q$ is simply $2$-connected and some ramification conditions, Hida obtained a proof of base change from $\Q$ to $F$ subject to a conjecture on permutation representations (cf. [H]).\\
\newline
\bf{Acknowledgments}: \rm
The first version of this paper was written during November 2009, when the Research Programme in Arithmetic Geometry at the Centre de Recerca Matematica (CRM) was taking place. I benefited from very useful conversations with several visitors, in particular A. Pacetti and R. Ramakrishna, so I want to thank them, and also the other visitors and  co-organizers of the Research Programme, and the staff of CRM. I also want to thank D. Geraghty for some comments on this new version. The referee's comments and remarks were also fundamental to improve the exposition, so it is a pleasure to thank him too.
\\
\newline
We finish this section with some definitions and notations.\\

\bf{Notation}: \rm In this paper, $F$ will always denote a totally real Galois number field. \\
For every number field $K$, we will denote by $G_K$ the absolute Galois group of $K$. \\
We will denote by $\chi$ the $p$-adic or mod $p$ cyclotomic character. The value of $p$, and whether it is the $p$-adic or the mod $p$ character, will always be clear from the context.\\
We will denote by $\omega$ a Teichmuller lift of the mod $p$ cyclotomic character.\\
Given a Galois representation $\sigma$, we will denote by $\mathbb{P}(\sigma)$ its projectivization.\\

\bf{Definitions}: \rm Let $K$ be a number field. Let $\bar{\rho}_p$ be a two-dimensional, odd, representation of $G_K$ with values on a finite extension of $\F_p$. 
\begin{enumerate}  \item We say that the image of $\bar{\rho}_p$ is \it{large} \rm if $p \geq 5$ and the image contains $\SL(2, \F_p)$. In this case, it is easy to see that the image of $\mathbb{P}( \bar{\rho}_p )$ is isomorphic to one of the following two groups: $\PSL(2, \F_{p^r})$ or $\PGL(2, \F_{p^r})$, for some $r$. Since $p \geq 5$, this implies in particular that large images are non-solvable.
	\item We say that the image of $\bar{\rho}_p$ is \it{dihedral} \rm when the image of $\mathbb{P}( \bar{\rho}_p )$ is a dihedral group. 
	\item We say that the image of $\bar{\rho}_p$ is \it{bad-dihedral} \rm when it is dihedral, $p >2$, and the quadratic number field where the restriction of $\bar{\rho}_p$ becomes reducible is $K(\sqrt{\pm p})$, where the sign is $(-1)^{(p-1)/2}$.
	\item We will say that $f$ is a \it{classical modular newform} \rm if it is a holomorphic newform of weight at least $2$.
\end{enumerate}

\section{General Description of the Proof}

We start with a given newform $f$ of odd level and $F$ a totally real Galois number field, as in Theorem \ref{teo:main}. We can assume that $f$ is not CM and that $\Gal(F/\Q)$ is non-solvable (otherwise, base change is known).  Let us make some remarks that apply to all steps of the proof. We will do a series of links as described in the introduction and we will manage to make all these links in such a way that
 in all the steps the residual representations will be irreducible, and their images will not be bad-dihedral.  We will also ensure that a similar restriction holds for the restrictions to $G_F$ of the residual images. Whenever applying the M.L.T. in [KiBT] and [Ki2], we will furthermore be in a situation where the residual projective image will contain some $\PSL(2,p^s)$ and elements of a prime order bigger than $5$ (thus these linear groups will be simple groups, and not $A_5$), and we will also see that the same holds for the restriction of the residual representations to $G_F$. In short: in most steps of our proofs (more precisely: as long as the good-dihedral prime is in the level, see the following paragraph for more details), we will have at the cha\^{i}nons a non-solvable residual image, even after restriction to $G_F$. \\
In order to get such a control on residual images, we will first introduce through level-raising a ``good-dihedral prime" $q$ in the level as in [KW1]. This is a technique created in the proof of Serre's conjecture given in [KW1] precisely to guarantee that residual images are non-solvable as long as the good-dihedral prime is not removed from the level. The next step is to apply the ``ramification swapping strategy" created in [D], in conjunction with suitable M.L.T. and the main result from [BLGG], to reduce to a situation where all primes in the level are split in $F$. After this step, all the characteristics where the congruences $f_i \equiv f_{i+1}$ hold will be primes that are split in $F/\Q$, in particular, all auxiliary primes introduced in the proof will be required to be split in $F/\Q$. \\
 We will then proceed to perform the killing ramification at primes in the level relying mainly on the M.L.T. in [KiFM]: to verify that the conditions to apply this result (to the restrictions to $G_F$ of the Galois representations that are the cha\^{i}nons of this chain, in {\bf reverse order}) hold, we will use ideas similar to those employed in [D] to control tame inertial weights of potentially crystalline representations under certain conditions (via results of Caruso, cf. [C]), together with a new, very useful, trick, that we will call the ``odd weight trick". \\
The main innovation of this paper comes at the step where (sooner or later the time of the farewell should come...) we remove the good-dihedral prime from the level. {\bf Before} doing so, we introduce a Micro-Good-Dihedral (MGD) prime to the level. This will be a small prime (in fact we will take $p=7$ in this paper, so let us just call it $7$) such that, after showing that by some level-raising arguments we can introduce ``supercuspidal" ramification at $7$, we are reduced to consider modular forms of level divisible by $49$ with ramification at $7$ being given by a character of order $8$ of the unramified quadratic extension of $\Q_7$. Thus, $7$ will work as a MGD prime in the following sense: as long as we work in characteristics $p$ such that $p \neq 2,7$ and $7$ is a square mod $p$ (and we will do so at all steps that go after losing our big good-dihedral prime $q$, except for  one step in characteristic $11$, where we will show why everything is fine by explicit computations) the local information at $7$ will be enough to ensure that the residual representations being considered are irreducible and not bad-dihedral, even after restricting to $G_F$.\\
\newline
We divide the long chain in three parts: In {\it Fase Uno} we introduce the big good-dihedral prime $q$ and, after reducing (via swapping) to a situation where all primes in the level are split in $F$ we play our ``odd weight trick" and kill ramification at all primes in the level. At the end we are reduced to the case of a newform with ``weight $k$, level $q^2$, good dihedral at $q$" with $k < q$. This newform is supercuspidal locally at $q$. In {\it Fase Dos} we play our level-raising trick to introduce the MGD prime $7$ in the level. After this, we play again the ``odd weight trick" (by introducing some nebentypus at $11$) and we are ready to kill ramification at $q$ (farewell, big good-dihedral prime!). When ending this second \it{fase} \rm we are in the case of a newform of level $7^2 \cdot 11$, with nebentypus at $11$ and some odd weight $k$. This newform is principal series locally at $11$, more precisely, the corresponding Galois representations have inertia Weil-Deligne parameter locally at $11$ equal to $(\psi \oplus 1, N = 0)$, where $\psi$ is the quadratic character corresponding to $\Q(\sqrt{-11})$, and it is supercuspidal locally at $7$. Finally, in {\it Fase Final}, we make the final moves so that our chain connects the original modular form with some newform $f_j$ of level $49$ and weight $k \leq 12$, trivial nebentypus (thus even $k$), and supercuspidal at $7$. For all such newforms (there are just two conjugacy classes of such newforms, one twisted of the other, in each of $k=4,6,8,10,12$, and none in $k=2$) we observe that a suitable member of the conjugacy class is ordinary at $3$ and has residual image in $\GL(2,\F_3)$, the image is known to be irreducible and not bad-dihedral even after restriction to $G_F$ (because of the local information at $7$). Thus, as in Wiles' first paper on modularity of elliptic curves (applying results of Langlands and Tunnell) we have mod $3$ modularity for the restriction to $G_F$ of these forms and applying  a M.L.T. for residually irreducible ordinary representations of Skinner and Wiles (cf. [SW]) we conclude that they can be lifted to Hilbert modular forms on $F$.\\
This completes the chain. In the next three sections we will go through the three \it{fases} \rm in full detail.\\

\section{{\it Fase Uno}}

Recall that all primes where we will build the cha\^{i}nons of our chain, except for the primes in the level of $f$, can be taken, and so will they be, to be split in $F/\Q$ (thanks in particular to the assumptions in Theorem \ref{teo:main}). Thus, this restriction will apply to all the auxiliary primes in the following construction.\\  

Since there will be several auxiliary primes in our proof, let us name them all right now, in particular to know their relative sizes. Let us call $p_i$, $i=1,..., w$ the prime factors of the level $N$ of $f$, and $k$ its weight. Recall that $N$ is odd, and that $f$ does not have CM.\\

 We will need auxiliary primes $b_i$, for $i=1,...., w$; $r_0, m$ and $r_1$, all split in $F$,  satisfying: $$ 2   \cdot \max \{ N,k \} < 2 \cdot b_i < r_0 < m < r_1 \qquad \; \quad (*)$$ for every $ 1 \leq i \leq w$. These are the auxiliary primes that we will use later. Let us take $B$ a constant bigger than $r_1$ (and always make sure in particular that $B$ is bigger than $7$).\\
\newline
Since we will apply the M.L.T. in [KiFM] to the restrictions to $G_F$ of our Galois representations at many steps, let us recall the statement of this Theorem. Since we will apply it to Galois representations that are base change of representations attached to classical modular forms, recall that these representations are known to be always potentially semistable (equivalently: de Rham) locally at $p$ and with different Hodge-Tate weights $\{0, k-1 \}$. Together with a condition on the size of the residual image, let us stress that
there is a {\bf technical condition} required for this Theorem to hold. 

\begin{teo}
\label{teo:superKisin} (Kisin) Let $F$ be a totally real Galois number field and $p$ be an odd prime.
Let $\rho$ be a representation of $G_F$ with values on a finite extension of ${\Q}_p$ that is $2$-dimensional, continuous, odd, absolutely irreducible and ramified at finitely many primes, with $p$ a prime that is split in $F$. Assume also that  the representation is, at all places $v \mid p$, de Rham of parallel Hodge-Tate weights $\{0, k-1 \}, \; k \geq 2$, and that the  residual representation $\bar{\rho}$ is absolutely irreducible even when restricted to the absolute Galois group of $F(\zeta_p)$. Let us also assume that the following {\bf technical condition} is satisfied: \\
         $\bar{\rho}|_{D_v}$ is not  isomorphic, for any $v \mid p$, to a twist of:
\[ \left( \begin{array}{ccc}
\chi & * \\
0 & 1   \end{array} \right)\] \\
         Then, if $\bar{\rho}$ is modular, $\rho$ is also modular.
\end{teo}

\bf{Remark}: \rm As recorded in a note added in proof to the published version of [KiFM], the condition that the representation should become semistable over an abelian extension of $\Q_p$ that appears in [KiFM] can now be removed due to recent work of Colmez.\\
\newline
Since whenever we apply this M.L.T. the representation will be base changed from $\Q$ and the characteristic $p$ will be split in $F$, it will suffice to check the last technical condition of the Theorem over $\Q$: whenever we are in a cha\^{i}non of our chain linking $f_i$ to $f_{i+1}$ and we are willing to propagate modularity over $F$ from $f_{i+1}$ to $f_{i}$ via Theorem \ref{teo:superKisin}, we will have to check that this local condition on the residual representation holds. \\
As for the condition on the residual image of the restriction to $G_{F(\zeta_p)}$, notice that whenever we can show that the restriction to $G_F$ has large image, since large implies non-solvable, this will be enough to see that this condition is satisfied.\\
\newline
We begin by changing to a weight $2$ situation, since this will be required in order to introduce a good-dihedral prime $q$ as in [KW1]. Given the family of Galois representations attached to $f$, since $f$ does not have CM we know that for almost every prime the residual image is large (due to Ribet's Theorem in [R1]) so we choose a characteristic $r_0$ where the residual image is large, split in $F$, and as in (*). We consider the residual representation, and take a minimal weight $2$ lift (as defined in [KW1]), corresponding to a weight $2$ modular form $f_2$ (if the weight of the given $f$ is $2$, this step is not taking place, we just ignore it and thus $f_2 = f$ in this case).  Thus, the newform $f_2$ has weight $2$ and its level is $N' \cdot r_0$, where $N'$ divides $N$. To ease notation, we will assume that $N' = N$.\\
\newline
Theorem \ref{teo:superKisin} will ensure that, when restricting to $G_F$, modularity propagates in reverse order (i.e., from $f_2$ to $f$). In order to check that the technical condition is satisfied just observe that $r_0$ is not in the level of $f$ and is bigger than twice its weight, thus due to Fontaine-Laffaille theory the residual tame inertia weights are equal to the Hodge-Tate weights $\{0, k-1 \}$. On the other hand, the residual image is large even after restriction to $G_F$ because of the following Lemma.

\begin{lema}
\label{teo:lema1} Let $p \geq 5$ be a prime and $\bar{\rho}$ a two-dimensional, odd, representation of $G_\Q$ with values on a finite extension of $\F_p$ and large image. Let $F$ be a totally real Galois number field. Then the image of the restriction of $\bar{\rho}$ to $G_F$ is also large, i.e., it contains the non-solvable group $\SL(2, \F_p)$.
\end{lema}

\bf{Proof}: \rm
Consider $\mathbb{P}( \bar{\rho} )$. We know that its image is of the form $\PSL(2, \F_{p^r})$ or $\PGL(2, \F_{p^r})$, for some $r$. Since $p>3$, the group $\PSL(2, \F_{p^r})$ is simple. If we consider the restriction of $\mathbb{P}( \bar{\rho} )$ to $G_F$, its image will be a normal subgroup of the image of $\mathbb{P}( \bar{\rho} )$, thus it either will be trivial or will contain $\PSL(2, \F_{p^r})$. \\
Since $\bar{\rho}$ is odd and $F$ is totally real, we know that the restriction to $G_F$ of $\mathbb{P}( \bar{\rho} )$ has non-trivial image, because the image of complex conjugation gives a matrix with eigenvalues $1$ and $-1$, thus non-trivial even modulo scalar matrices. \\
Therefore, we conclude that the image of $\mathbb{P}( \bar{\rho} )$ must contain $\PSL(2, \F_{p^r})$, and this proves the Lemma.
\qed\\
\newline
Now we want to add the good dihedral prime to the system of Galois representations attached to $f_2$. More precisely, we need to find
 two primes $t$ and $q$ greater than $B$ such that, as in [KW1] and [D], modulo $t$ we can do level raising to introduce the extra ramification at $q$ so that in the next steps the Galois representations in characteristics smaller than $B$ will have the good dihedral prime $q$ in their ramification set. Since we will require that these two primes are also split in $F$, let us give in the following Lemma the definition and proof of existence of these two primes. 

\begin{lema}
\label{teo:lema2} Let $F$ be a totally real number field. Let $\{ \rho_\ell \}$ be a compatible system of Galois representations attached to a classical newform $f$ of weight $2$ and level $N$, such that $f$ does not have CM. Let $B$ be a constant greater than $N$ and $7$. Then:\\
 There is a prime $t > B$ such that $t \equiv 1 \; \mod{4}$, t splits in $F$, and the image of $\bar{\rho}_t$ is exactly $\GL(2, \F_t)$.\\
 Furthermore, there is a prime $q$ satisfying the following conditions:\\
 
\begin{enumerate}
	\item $q \equiv -1 \; \mod{t} $
	\item The image of $\bar{\rho}_t(\Frob \; q)$ has eigenvalues $1$ and $-1$
	\item $q \equiv 1 \; \mod{8} $
	\item $q \equiv 1 \; \mod{p} $ for every $p \leq B$
	\item $q$ splits in $F$
\end{enumerate}
\end{lema}

\bf{Proof}: \rm
A similar result is proved in [KW1] and [D], without the condition that $t$ and $q$ be split in $F$. Let us explain why the result is true with this extra condition.\\
The existence of $t$ follows from the result of Ribet in [R1], which implies that for almost every prime the residual images of the modular compatible system will be large. Then it suffices to take $t$ sufficiently large and split in the compositum of $\Q(i)$, $F$, and the field of coefficients $\Q_f$ of $f$. Observe that the determinant of $\bar{\rho}_t$ is $\chi \cdot \varepsilon$, where $\varepsilon$ is some Dirichlet character unramified at $t$.  \\
The existence of $q$ was proved using Cebotarev density Theorem in [KW1] and [D] without the requirement that $q$ splits in $F$, thus we can deduce that the Lemma is true using again Cebotarev if we see that the extra condition (5) is compatible with the other conditions on $q$. This is immediate for conditions (3) and (4), but not so obvious for conditions (1) and (2).\\
Since $t$ is split in $F$, $F$ is linearly disjoint from $\Q(\zeta_t)$, and this implies that condition (5) and (1) are compatible.\\
In [KW1] a prime $q$ satisfying (2) is obtained by taking $\Frob \; q$ in the same conjugacy class of complex conjugation. Since $F$ is totally real, the same argument works if we work with $G_F$, thus proving that conditions (5) and (2) are compatible. \\
Therefore, by Cebotarev's density Theorem, we conclude that there exists a prime $q$ satisfying all the conditions in the statement of the Lemma.
\qed\\
 \newline
 
 We apply this Lemma to the system attached to $f_2$, taking
 $B$ to be the constant defined at the beginning of this section. \\
Then, following [KW1], we know that there is a congruence mod $t$ between $f_2$ and a newform $f_3$ of weight $2$ such that $q^2$ divides the level of $f_3$ and ramification at $q$ of the attached Galois representations is given by a character of the quadratic unramified extension of $\Q_q$ of prime order $t$. \\
\bf{Remark}:  \rm In fact, [KW1] constructs this non-minimal lift using potential modularity, and this can be done because the conditions to apply M.L.T. are satisfied, now in our case since by assumption this mod $t$ representation of $G_\Q$ is modular, it clearly follows that the non-minimal lift is also modular, thus attached to some newform $f_3$ of weight $2$ (alternatively, we could have used the theory of congruences between modular forms and raising the level). \\
\newline
For the restrictions to $G_F$, we will need to go in reverse order: since the residual image is large even after restriction to $G_F$ because of Lemma \ref{teo:lema1}, and both $t$-adic representations are Barsotti-Tate at $t$, over $F$ modularity propagates from $f_3$ to $f_2$ due to the M.L.T. in [KiBT] (note that, as observed by Kisin in the first section of [Ki2], in cases where $t$ is totally split in $F$ this applies without having to check that ``pot. ordinary" goes to ``pot. ordinary").\\
\newline
\bf{Very Important Remark}: 
\rm From now on, as long as we work in characteristic $p \leq B$, we know that all representations in our chain will be residually irreducible and not dihedral (because of the local information at the good dihedral prime $q$), moreover they will have projective image containing some simple group $\PSL(2, \F)$ (not isomorphic to $A_5$). This follows from the properties of a good dihedral prime and a study of images in this case (see [KW1], Lemma 6.3). \\
\newline
 Let us show that we can also control the restriction to $G_F$ of the images:
\begin{lemma} \label{teo:lema1.5}
Let $p$ be any prime smaller than the bound $B$ introduced above. Let $\bar{\rho}_p$ be a $2$-dimensional odd representation of $G_\Q$ with values on a finite extension of $\F_p$. Let $q > B$ be a good dihedral prime for the representation $\bar{\rho}_p$, and let $F$ be a totally real number field such that $q$ is split in $F$. \\
Then the image of $\mathbb{P}(\bar{\rho}_p)$ restricted to $G_F$ contains a simple group of the form $\PSL(2, \F_{p^r})$.
\end{lemma}
\bf{Proof}: \rm
Since $q$ is split in $F$, when restricting to $G_F$ the image of the
 decomposition group at $q$ (this restriction is dihedral by assumption)
 forces the residual image to stay irreducible. Thus, since (as we have observed above) the good-dihedral prime forces the image over $\Q$ of the projectivization to be an almost simple linear group and $F$ is Galois, we see that even after restricting to $G_F$ the image will contain a simple group $\PSL(2, \F_{p^r})$. 
 \qed\\
\newline
Thus, as long as we work in characteristic $p < B$ the results in [Ki2], [KiBT] and Theorem \ref{teo:superKisin} will allow us to show that modularity over $F$ propagates in reverse order through the chain, if we can show that the technical condition in Theorem \ref{teo:superKisin} is satisfied and if we can deal with the primes in the level that are not split in $F$. In other words, we don't have to worry about the size of the restriction to $G_F$ of the residual image thanks to Lemma \ref{teo:lema1.5}. Our next step is to reduce to a situation where all primes in the level are split in $F$.\\

\bf{Ramification Swapping}: \rm \\ 
Before starting, let us make a {\bf General Remark} that applies to all the steps of the proof, not only to this section:  as in [KW1] and [D], we always assume that residual representations in characteristic $p$ have Serre's weight $k \leq p+1$, because it is well-known that by making a suitable twist one can reduce to this case, and twists preserve modularity.  \\
\newline
Swapping is the process used in [D], section 4, in order to transfer ramification from one set of primes to another. We have currently the set of primes: $p_1, ..., p_w$ in the level, together with $r_0$ and the good-dihedral prime $q$. 
To ease notation, we assume that none of the primes $p_i$ is split in $F$. Therefore, because of the assumptions of Theorem \ref{teo:main}, we can suppose that they are all greater than $5$. The auxiliary primes $b_1, ....,b_w$ are chosen to be split in $F$ and satisfying the inequalities (*) described at the beginning of this section. In particular, they are larger than all the $p_i$ but smaller than the bound $B$. Thus, in all what follows we know that residual images will be large even when restricting to $G_F$, due to Lemma \ref{teo:lema1.5}.\\
Starting with the newform $f_3$, whose weight is $2$, let us recall the process of ramification swapping from [D]: we move to characteristic $p_1$, reduce mod $p_1$, and take a minimal lift. In general, this lift will correspond to a newform $f_4$ of weight $k>2$. What is more important: the prime $p_1$ is not in the level of $f_4$. Modularity is preserved, for the restrictions to $G_F$, from the Galois representation attached to $f_4$ to the one attached to $f_3$, as follows from the M.L.T. of [G], [KiBT], and [BLGG]. More precisely, we show this by dividing in two cases:

\begin{enumerate}

	\item The Galois representation attached to the weight $2$ newform $f_3$ is potentially Barsotti-Tate at $p_1$: In this case, the M.L.T. in [KiBT] can be applied, but it requires the construction of ordinary modular lifts. This was accomplished in [BLGG], Theorem 6.1.11, where they conclude that if the prime $p_1$ is not  split in $F$ (in the split case the result of [KiBT] applies automatically) then  the M.L.T. of Kisin can be applied over $F$ if the following condition is satisfied: $[F(\zeta_{p_1}) : F] > 4$. Since we have $p_1 >5$, this is satisfied if the prime is unramified in $F$, but may fail for ramified primes.  What we do to remedy this situation is to apply  
	solvable base change (cf. [L] and [AC]): 
	if we consider a subfield $F'$ of $F$ such that the Galois group $Gal(F/F')$ is solvable, and a representation of $G_{F'}$, then it is known that the modularity of such a representation is equivalent to the modularity of its restriction to $G_F$, assuming that this restriction is irreducible. Let us see that this implies that we can replace $F$ by $F'$ at this step of the chain:\\
	We start with the assumption that $\rho_{f_4, p_1}|_{G_F}$ is modular. Then, by solvable base change (``going down"), $\rho_{f_4, p_1}|_{G_{F'}}$ is also modular. We have a congruence modulo $p_1$ between $f_3$ and $f_4$, and suppose that for the restrictions to $G_{F'}$ of the Galois representations attached to $f_3$ and $f_4$ modularity propagates well from $f_4$ to $f_3$. Then we conclude that $\rho_{f_3, p_1}|_{G_{F'}}$ is also modular. Finally, by a second application of solvable base change (``going up"), we conclude from this that $\rho_{f_3, p_1}|_{G_F}$ is modular.\\
	In other words: we can replace $F$ by any subfield $F'$ of it such that $Gal(F/F')$ is solvable and if we can show that applying a suitable M.L.T. modularity propagates well for the restrictions to $G_{F'}$, then the same holds for the restrictions to $G_F$.\\
	So, let us choose the right subfield of $F$. Let $v$ be a prime of $F$ dividing $p_1$, and let $D$ be the decomposition group of $v$ in $Gal(F/\Q)$. Let $F^{(v)}$ be the fixed field of $D$. Then, since $D$ is solvable, we can take this as a candidate subfield of $F$ to apply the above strategy.  The field $F^{(v)}$, by construction, has a place dividing $p_1$ that is split over $\Q$, thus $F^{(v)}$ is linearly disjoint from $\Q(\zeta_{p_1})$, and therefore the condition $[F^{(v)}(\zeta_{p_1}) : F^{(v)}] > 4$  is clearly satisfied (recall that $p_1 > 5$).  Recall that we know that the residual image restricted to $G_F$ is large, and a fortiori the same is true over $G_{F^{(v)}}$. Therefore, we can apply  Theorem 6.1.11 in [BLGG] and conclude that for the restrictions to  $G_{F^{(v)}}$, modularity propagates well from $f_4$ to $f_3$, and by the above ``solvable base change trick" that the same is true also over $G_F$.  
	
	\item The Galois representation attached to the weight $2$ newform $f_3$ is potentially semistable at $p_1$: this case corresponds to the local component at $p_1$ of $f_3$ being Steinberg or twist of Steinberg. In particular, by taking a suitable twist, we can assume that in this case the Galois representation is semistable at $p_1$. Then, the M.L.T. of Geraghty (cf. [G], Theorem 5.4.2) applies because this $p_1$-adic Galois representation is ordinary, and the one attached to $f_4$, which in this case will be crystalline of weight $2$ or $p_1 + 1$ (depending on the Serre's weight of the mod $p_1$ representation), since it is forced by the congruence with $f_3$ to be residually ordinary it is also known to be ordinary (this is well-known for weight $2$ newforms, and follows from the results in [BLZ] in the case of weight $p_1+1$) and, a fortiori, their restrictions to $G_F$ are also ordinary at all places dividing $p_1$. In order to apply this M.L.T., there are two other technical conditions that need to be checked, which appear as conditions (4) and (5) in Geraghty's Theorem (cf. [G], Theorem 5.4.2). 
	\\ Condition (5) requires the image of the residual representation restricted to $G_{F(\zeta_{p_1})}$ to be ``big", a technical notion which, for the case of $2$-dimensional Galois representations,  is known to hold whenever this image is large (see for example [BLGG], section 4). We know that our residual representation has large image when restricted to $G_F$, and therefore it also has large image when restricted to $G_{F(\zeta_{p_1})}$ (because $F(\zeta_{p_1})$ is a cyclic extension of $F$), thus condition (5) of Geraghty's Theorem is satisfied. 
	\\ Condition (4) is the requirement that the extension of $F$ fixed by the kernel of the adjoint of
the residual representation does not contain the $p_1$-th roots of
unity. Since  we know that the residual image of the restriction to $G_F$ is large, it is easy to 
see (this argument appears in the proof of Theorem 6.1.9 in [BLGG])
that this condition  is implied by the following one:  $[ F(\zeta_{p_1}) :
F  ] > 4$. Therefore, we can just apply the ``solvable base change trick" as we did in the previous case to reduce to a situation where
this condition is satisfied (it is obvious that in this change of field we are preserving largeness of the residual image, and also that the $p_1$-adic representations remain ordinary because they are restrictions of ordinary representations of $G_\Q$). 
We conclude that the M.L.T. of Geraghty
applies in this case and allows us to propagate modularity, for the restrictions to
$G_F$, from $f_4$ to $f_3$.\\
{\bf Remark}: In fact, a stronger version of the theorem of Geraghty is
given in [BLGGT], Theorem 2.3.1 (it is written for imaginary CM fields, but a
standard and very easy argument using quadratic base change, appearing
for example in [G], shows that a similar result also holds over
totally real fields). In this version, results of Thorne have been
incorporated and thanks to them condition (4) in Geraghty's theorem
disappears and condition (5) is replaced by a condition which is known
to be satisfied if the prime is greater than $6$ and  the restriction of the residual representation to $G_{F(\zeta_{p_1})}$ is irreducible. Again,
since the image of the restriction to $G_F$ of our residual
representation is large, this is clearly satisfied in our case.

\end{enumerate}

With the newform $f_4$, of weight $ 2 \leq k \leq p_1 + 1$ (and level prime to $p_1$) we move to characteristic $b_1$, reduce modulo $b_1$, and take a minimal weight $2$ lift corresponding to a modular form $f_5$. Observe that by construction the prime $b_1$ is not in the level of $f_4$, hence the $b_1$-adic representation attached to $f_4$ is crystalline. By the inequalities (*), the weight of $f_4$ is much smaller than $b_1$, thus because of  Fontaine-Laffaille theory and the fact that $b_1$ is split in $F$ we see that the technical condition in Theorem \ref{teo:superKisin} is satisfied, thus modularity for the restrictions to $G_F$ can be propagated from $f_5$ to $f_4$. At this step we are assuming that the weight of $f_4$ is greater than $2$: if this is not the case, we just take $f_5 = f_4$ and there is no need of adding ramification at $b_1$. For simplicity, let us just assume from now on that $b_1$ is in the level of $f_5$ in any case, since this will not affect the rest of the proof (i.e., whether or not $b_1$ is truly on the level we proceed in the same way and the outcome is the same). \\
What remains is just an iteration of the above procedure: notice that $f_5$ is again of weight $2$, it has $b_1$ in its level, and the  level is prime to $p_1$. Thus, we repeat the procedure, for every $1 < i \leq w$, moving first to a characteristic $p_i$ with a weight $2$ family, reducing and taking a minimal lift (this kills the ramification at $p_i$), and then to  characteristic $b_i$, reducing and taking a minimal weight $2$ lift (thus introducing ramification at $b_i$). As we have indicated, via M.L.T. of Geraghty, Kisin, and the result in [BLGG], modularity of the restrictions to $G_F$ is preserved through the whole process.\\
The process concludes when one makes the last minimal weight $2$ lift, in characteristic $b_w$. We end up with a newform $g$ of weight $2$ whose level does not contain any of the primes $p_i$, and in their place we have the primes $b_i$ which are split in $F$.\\
In fact, as we explained in the first iteration, maybe not all the $b_i$ appear in the level of $g$, but for simplicity we will just act as if they do: in particular, in the ``iterated killing ramification" step we will move to each characteristic $b_i$ even if it is not in the level (this is still useful because by twisting we may change the weight in each of these characteristics).\\
 We have reduced the proof of our main Theorem to a situation where all primes in the level are split in $F$. Recall also that from now on, all the residual characteristics that will appear in the chain are going to be split in $F$.\\
 \newline
\newline
Now we want to manipulate a bit the nebentypus of $g$. The following step is meant to reduce to a situation where, for every prime $b_i$ or $r_0$ in the level such that the nebentypus ramify at it, we have that the corresponding abelian extension of $\Q$ (contained in $\Q(\zeta_{b_i})$) has odd degree, and is thus real. To achieve this, we simply move to characteristic $2$ (recall that by assumption $2$ is not in the level of $f$, thus it is also not in the level of any newform in our chain). Since $g$ has weight $2$, we can reduce mod $2$ and take a minimal lift (as in [KW1]), which will correspond to another weight $2$ modular form $g_2$. Since the lift is minimal (locally at every prime) we are reduced to a situation where the $b_i$-part and the $r_0$-part of the nebentypus is ``real" (because it has odd degree) for any prime dividing its conductor. This is because a character with values on a finite extension of $\F_2$ must have odd order. Observe that we want that, on the restriction to $G_F$, modularity can be propagated from $g_2$ to $g$, and this follows from the main result in [Ki2] since $g$ is Barsotti-Tate at $2$ and $2$ is split in $F$ (and the residual image is, as follows from \ref{teo:lema1.5}, non-solvable even when restricted to $G_F$).\\
\newline
Let us introduce an auxiliary prime $m$ whose role will be to produce ``odd Serre's weights" for the residual representations to be dealt with in the ``iterated killing ramification" step, where we kill ramification at all primes $b_i$ and $r_0$ in the level. The prime $m$ is a prime as in (*): it is smaller than $B$, split in $F$, it is congruent to $3$ mod $4$, and bigger than twice all of the $b_i$ and bigger than $r_0$. We move to characteristic $m$ and here we reduce the weight $2$ Galois representation attached to $g_2$ and take a lift given by another newform $g_3$ which has nebentypus $\psi$ of order $2$ corresponding to the quadratic extension $K= \Q(\sqrt{-m})$, and a suitable weight $k > m^2 - 1$. It is obvious that we can find a $k > m^2-1 $ such that the congruence mod $m$: $\psi \cdot \chi^{k-1} \equiv \chi$ holds, and having this, it is a result of Khare (cf. [Kh], Theorems 1 and 2) that a congruence between $g_2$ and a newform $g_3$ whose weight and type at $m$ are $(k, \psi \oplus 1, N=0)$ does hold. Since $K$ is imaginary it is clear that $k$ will be odd (because modular Galois representations are odd).  Modularity propagates well on the restriction to $G_F$ from $g_3$ to $g_2$ because of [KiBT]. \\
\newline
Now we simply kill ramification at the primes $b_i$, and $r_0$, as in the ``iterated killing ramification" (I.K.R.) step in [D]: by switching to each of them, reducing mod $b_i$, and taking a modular minimal lift (recall that by suitable twisting we can assume that it will have weight $k_i \leq b_i + 1$), then moving to the next, and so on.  Since the nebentypus at each $b_i$ and at $r_0$ is at most given by a real abelian extension (we have managed to reduce to such a case), it is clear that we start with $g_3$ of odd weight and the Serre's weight $k_i$ mod $b_i$ will  also be odd, and this is enough to see that the technical condition in Theorem \ref{teo:superKisin} holds true. As usual we want to propagate modularity for the restrictions to $G_F$ of these representations, in reverse order, and since we know that residual images are large even when restricting to $G_F$ (by Lemma \ref{teo:lema1.5}) it is enough to verify this condition. During all the I.K.R.  the residual representations will have odd Serre's weight and thus Theorem \ref{teo:superKisin} can be applied. When we finish, we end with a newform $g_s$ of level $m \cdot q^2$ and odd weight $k_s \leq b'$ ($b'$ being the smallest of the primes $b_i$ in the level, since we can perform I.K.R. with the primes taken in decreasing order), thus $m$ is bigger than twice this weight because of the inequalities (*). We now move back to characteristic $m$ and we consider the residual representation of $g_s$. This one will have even Serre's weight, but since ramification at $m$ was just given by a character of order $e= 2$ and we are reaching characteristic $m$ with a family of weight $k_s$ smaller than $m/2$,  using the results of Caruso in [C] as we did in [D] one can check that in this situation the technical condition in Theorem \ref{teo:superKisin} is satisfied.\\
 In fact, if we extract what was proved in [D] using the results in [C] (this is contained in the proof of Lemma 4.4 in [D]) we have the following:
 
 \begin{lema} \label{teo:caruso}
 Let $\rho_p$ be an odd, continuous, two-dimensional representation of $G_\Q$, with finite ramification set and values on a $p$-adic field. Suppose that locally at $p$ the representation is potentially crystalline, and that it becomes crystalline when restricted to a subfield of $\Q_p(\zeta_p)$ of even ramification degree $e$. Suppose that the Hodge-Tate weights of $\rho_p$ are $\{ 0, k-1 \}$ with $k>2$. Furthermore, assume that the following inequality is satisfied:
 $$ (k-1) \cdot e  <  p-1  $$
 Then, if we consider the restriction to the decomposition group at $p$ of the residual representation $\bar{\rho}_p$, the technical condition in Theorem \ref{teo:superKisin} is satisfied.
 \end{lema}
 
 Then, since we have $ e \cdot k_s = 2 \cdot k_s < m$ and $k_s>2$ (because $k_s$ is odd), thus in particular $ 2 \cdot (k_s -1) < m - 1$, we see from this Lemma that the M.L.T. of Kisin can be applied.
 This means that if we take a minimal modular lift corresponding to some newform $g_{s+1}$, modularity over $F$ propagates well from $g_{s+1}$ to $g_s$.\\
We have thus ended this {\it fase} with a newform $g_{s+1}$ of level $q^2$, trivial nebentypus, some even weight $k < B$, and good-dihedral at $q$.

\section{{\it Fase Dos}}
The purpose of this {\it fase} is to introduce some extra ramification at $7$ to the level as described in section 2. The type of ramification at $7$ will be similar to the one in $q$, except that it will correspond to a character of even order $8$. \\
To ease notation, let us rename $g_{s+1} := h_1$.\\
We begin, as we did in the previous section, by changing to a weight $2$ situation (we do this if and only if we have $k>2$). This is where our last auxiliary prime in the sequence (*) appears: it is a prime $r_1$ split in $F$, bigger than $m$ (thus bigger than the weight $k$ of $h_1$) and smaller than $B$. We move to characteristic $r_1$, consider the residual representation and take a modular weight $2$ minimal lift corresponding to the newform $h_2$ . Observe that since $k$ was even, the nebentypus at $r_1$ that we introduce here corresponds to a real abelian field. Residual images are large even after restricting to $G_F$  because of Lemma \ref{teo:lema1.5}, and Theorem \ref{teo:superKisin} (for the restriction to $G_F$, in reverse order) can be applied (this is due to the results of Fontaine-Laffaille, as in the case of characteristic $r_0$ in section 3). \\
The newform $h_2$ has weight $2$, and it has level $r_1 \cdot q^2$. At the prime $q$ the ramification is, as usual, supercuspidal, and at the prime $r_1$ the ramification, introduced with the weight $2$ lift, is given by the character $\omega^{k-2}$, i.e., the inertial Weil-Deligne parameter at $r_1$ for the Galois representations attached to $h_2$ is $(\omega^{k-2} \oplus 1, N = 0)$.\\
Now we move to characteristic $3$, reduce mod $3$, and take a weight $6$ modular minimal lift  of (some twist of) this residual representation corresponding to a newform $h_3$. Such a weight $6$ lift always exists (cf. [E]; see also [RS], Theorem 2.6). Observe that since $h_2$ was of weight $2$ and level prime to $3$, its $3$-adic representation is Barsotti-Tate, thus modularity of the restrictions to $G_F$ propagates from $h_3$ to $h_2$ by [KiBT].\\
The newform $h_3$ has weight $6$, and its level is again $r_1 \cdot q^2$, moreover, at both primes in the level the local Weil-Deligne inertial parameter is the same as the one of $h_2$.\\
Now we move to characteristic $7$, reduce mod $7$, and take a weight $2$ modular form $h_4$ lifting it: Theorem \ref{teo:superKisin} ensures that our chain works well at this step, for the restrictions to $G_F$, from $h_4$ to $h_3$, since $7>6$ and thus the results of Fontaine-Laffaille apply\footnote{in case the residual representation is, locally at $7$, reducible and decomposable, it will be (at least on inertia) isomorphic to the sum $\chi^{5} \oplus 1$ and this is a twist of $\chi \oplus 1$. In this case instead of applying Theorem \ref{teo:superKisin} we can apply the M.L.T. in [SW]. In fact, since $6<7$ the $7$-adic crystalline representation of Hodge-Tate weights $\{0, 5 \}$ on $h_3$ is known to be ordinary in this case, and the same is also known for the one corresponding to $h_4$, which is potentially Barsotti-Tate, in this residually ordinary case}.\\
The newform $h_4$ has ramification at $7$ given by the character $\omega^4$ of order $3$, i.e., for the Galois representations attached to $h_4$ the Weil-Deligne inertial parameter at $7$ is $(\omega^4 \oplus 1, N = 0)$. Now we move to characteristic $3$ as in the Sophie Germain trick in [D] (this works because the primes $3$ and $7$ are a pair of Sophie Germain primes): since $\omega^4$ has order $3$ the mod $3$ representation attached to $h_4$ will be either unramified or semistable at $7$. Moreover, in the first case using the fact  that $\rho_{h_4, 3}$ restricted to the decomposition group at $7$ is isomorphic to $\omega^4 \oplus 1$, and that the order $3$ character $\omega^4$ trivializes when reduced modulo $3$, we see that Ribet's sufficient condition for (semistable) raising-the-level at $7$ holds (cf. [R2]), in fact the image of $\Frob \; 7$ for the residual mod $3$ representation of $h_4$ has the eigenvalue $1$ with multiplicity $2$ and $7 \equiv 1 \pmod{3}$. \\
This means that in any case we have a modular weight $2$ lift corresponding to a newform $h_5$ with semistable ramification at $7$, and the chain works well at this step over $F$ due to [KiBT] since $h_4$ is Barsotti-Tate at $3$. \\
Now we use the fact that $7 + 1 = 8$ and we consider the mod $2$ representation attached to $h_5$. Because $h_5$ has semistable ramification at $7$, we are in a case where the results in [KW1] can be applied to produce a non-minimal weight $2$ modular lift of this mod $2$ representation (as usual, since in [KW1] they rely on potential modularity, thus on M.L.T., the lifts produced using their techniques, with the residual modularity assumption, are automatically modular), which is non-minimal only at $7$: we can ensure that this lift  has ramification at $7$ given by a character of order $8$ of the unramified quadratic extension of $\Q_7$ (we are taking $j=0$ and $i=6$ in the notation of [KW1]: this character of order $8$ is $\psi_2^6$, where $\psi_2$ denotes a fundamental character of level $2$ of the tame inertia group at $7$). Let us call $h_6$ the weight $2$ newform with level $49 \cdot q^2 \cdot r_1$ just produced. Observe that the attached residual Galois representations, in each characteristic $p \neq 2,7$, will have the MGD (Micro-Good-Dihedral) prime $7$ and in particular ramification at $7$ corresponding to a degree $8$ character of the unramified quadratic extension of $\Q_7$. Note that this character has order $4$ in the projectivization of the representations.\\

Remark: \bf{MGD prime: what is it and how does it help}: \rm a MGD prime is a prime $s$ in the level of a residual representation $\bar{\rho}$ in characteristic $p$ such that locally at $s$ it has the same local parameter as in the definition of good dihedral prime (i.e., induced from a character of an unramified quadratic extension of $\Q_s$), but without any further relation between $s$ on one hand and $p$ and the primes in the level on the other (this is the main difference with good dihedral primes). In particular, having a MGD prime in the level implies (with the same proof used for good dihedral primes) that the residual representation is irreducible. Furthermore, if it happens to be the case that the MGD prime $s$ is a square mod $p$, then (again, with the same proof used to control the images using good dihedral primes) the residual image can not be bad-dihedral.\\
\newline
Modularity of the restrictions to $G_F$ can be propagated from $h_6$ to $h_5$ using [Ki2] because $h_5$ is Barsotti-Tate at $2$ (and residual images are non-solvable due to the good-dihedral prime $q$, even after restriction to $G_F$, see Lemma \ref{teo:lema1.5}).
To finish this section, we will kill ramification at $r_1$ and then at $q$. In order to do so, we need again the trick of ``odd Serre's weights" in order to ensure that the chain propagates well modularity (as usual, over $F$ and in reverse order) via Theorem \ref{teo:superKisin}. Thus, we consider the auxiliary prime $11$ (also split in $F$ by assumption) and we reduce $h_6$ mod $11$ and take a modular odd weight lift $h_7$ of some weight $k > 11^2 -1$ and nebentypus given by the character of order $2$ ramifying only at $11$ (we rely again on [Kh]). The weight $k$ is odd because the field $\Q(\sqrt{-11})$ is imaginary. Since $11$ is not in the level of $h_6$ and $h_6$ is a weight $2$ newform, once again the results in [KiBT] allow to propagate modularity over $F$, from $h_7$ to $h_6$.\\
As we did in the previous section, since the weight is odd and the nebentypus at $r_1$ is real (and there's no nebentypus at $q$) we just move to characteristic $r_1$, reduce mod $r_1$ and take a minimal modular lift corresponding to a newform $h_8$, then the residual representation has odd Serre's weight and thus modularity can be propagated (over $F$, in reverse order) because of Theorem \ref{teo:superKisin}. \\
Then, we do the same in characteristic $q$: we take the modular form $h_8$, move to characteristic $q$, reduce mod $q$, then take a minimal modular lift corresponding to some newform $h_9$.    The technical condition in Theorem \ref{teo:superKisin} is satisfied once again because of the odd weight trick. But at this last step we have to be careful (for the first time!) with the residual image: since we are losing the good-dihedral prime $q$, we may have a small residual image. Here is where the MGD prime $7$ starts playing his role. Since $q \equiv 1 \pmod{8}$ and $q$ is a square mod $7$, then $7$ is a square mod $q$. Then, the usual arguments with good-dihedral primes (cf. [KW1]) imply that because of the decomposition group at $7$ having dihedral image this mod $q$ representation must be irreducible (it is so locally at $7$) and, if it is dihedral, it is not bad-dihedral (because $7$ is a square mod $q$). Moreover, since $7$ is split in $F$, then we also have that after restriction to $G_F$ the projective representation contains in its image a group that is dihedral of order $8$ (namely, the image of the decomposition group at $7$). In this situation, the following Lemma shows that the size of the residual image is good enough to apply Kisin's Theorem \ref{teo:superKisin}. 

\begin{lema} \label{teo:lema3} Let $F$ be a totally real Galois number field. If $p \geq 5$ is a prime and $\bar{\rho}_p$ is a two-dimensional, odd, representation of $G_\Q$ with values on a finite extension of $\F_p$ such that $\bar{\rho}_p$ is irreducible and its image is not bad-dihedral, and such that $\mathbb{P}(\bar{\rho}_p)$ restricted to $G_F$ contains a dihedral subgroup of order $8$, then the restriction of $\bar{\rho}_p$ to the absolute Galois group of $F(\zeta_p)$ is absolutely irreducible.
\end{lema}
\bf{Proof}: \rm
Using Dickson's classification of maximal subgroups of $\PGL(2, \F_{p^r})$ we see that the assumptions on $\bar{\rho}_p$ imply that its image must be of one of the following types:\\
\begin{itemize}
	\item (i) Large, i.e., containing $\SL(2, \F_p)$
	\item (ii) projectively isomorphic to $S_4$ or $A_5$
	\item (iii) dihedral, but not bad-dihedral.
\end{itemize}
In case (i), we apply Lemma \ref{teo:lema1} and conclude that the image of the restriction of $\bar{\rho}_p$ to $G_F$ is also large, thus containing a non-solvable group. In particular, the restriction to $F(\zeta_p)$, a cyclic extension of $F$, can not be reducible.\\
In case (ii), we use the assumption that the restriction of the projective image to $G_F$ contains a dihedral subgroup of order $8$, and the facts that $A_5$ is simple and $S_4$ does not contain a normal subgroup of order $8$ to deduce that the projective image does not change when restricting to $G_F$. Thus, being the projective image restricted to $G_F$ as in (ii) it is clear again that the restriction to the cyclic extension $F(\zeta_p)$ of $F$ can not be reducible.\\
In case (iii), let us call $K$ the quadratic number field such that the representation $\bar{\rho}_p$ is induced from a character of $G_K$. By assumption the image of the restriction to $G_F$ contains a dihedral group, thus we see that the restriction to $G_F$ is dihedral. This implies that $F$ is linearly disjoint from $K$.
Also, $K$ is not the quadratic number field contained in $\Q(\zeta_p)$ (this is the assumption of the representation being not bad-dihedral). \\
We conclude that $K$ is linearly disjoint from $F(\zeta_p)$, thus $\bar{\rho}_p$ restricted to $F(\zeta_p)$ is absolutely irreducible.
\qed\\
\newline

Therefore, having checked that Theorem \ref{teo:superKisin} applies (because of the Lemma above and the odd weight trick), we see that modularity of the restriction to $G_F$ propagates well from $h_9$ to $h_8$.\\ 
We end up with a newform $h_9$ of some odd weight $k$, level $7^2 \cdot 11$, quadratic nebentypus at $11$ and such that $7$ is a MGD prime for it. We should try (but we will not always be able to do so), in our next moves, to work in characteristics $p$ such that $7$ is a square mod $p$, since MGD primes work better there, i.e, they allow to conclude that the residual image is not bad-dihedral (as we did above in characteristic $q$).

\section{{\it Fase Final}}
We consider the mod $3$ representation of $h_9$. Since $3$ is not in the level and because of the nebentypus at $11$ it will have odd Serre's weight (thus it will have, up to twist, $k=3$). Since $7$ is a square mod $3$, we see from the MGD prime $7$ that this residual representation is irreducible and it is not bad-dihedral. We take a minimal modular lift, it corresponds to a newform $h_{10}$ of weight $3$  and level $7^2 \cdot 11$, quadratic nebentypus at $11$ and such that $7$ is a MGD prime for it. Because the residual Serre's weight is odd we can apply Theorem \ref{teo:superKisin} to propagate modularity over $F$ from $h_{10}$ to $h_9$. Just observe that because of the MGD prime $7$ and the fact that $7$ is split in $F$ the restriction to $G_F$ of this mod $3$ representation can not be reducible (it is irreducible locally at $7$) and because $7$ is a square mod $3$ it will stay irreducible if we restrict to $F(\sqrt{-3})$, as required in Theorem \ref{teo:superKisin}.\\
Now we consider the mod $11$ representation attached to $h_{10}$: it will be irreducible because it is so locally at $7$, but unfortunately $7$ is not a square mod $11$, thus a priori it could be bad-dihedral.
To check that it is not bad-dihedral, we do some computations. Observe that we are dealing with a mod $11$ modular representation of level $49$, some even weight that can be taken (by twisting) to be $k \leq 12$, and such that it is supercuspidal at $7$. We check in W. Stein's tables that, except for $k=2$ where there is not any such newform, there are two conjugacy classes of newforms in each of the other spaces satisfying (in fact: that may satisfy) these conditions, and they are twists of each other. Suppose that the mod $11$ representation of $h_{10}$ is bad-dihedral. Then, it is well-known (and can be easily proved by looking at the action of the inertia group at $11$ and using the definition of Serre's weight) that this can only happen if $11 = 2k - 3$ or $11 = 2k -1$ (cf. [KW1], Lemma 6.2). Since $k$ is even, the only possibility is thus $k=6$. Reducing eigenvalues mod $11$ we easily check that for the couple of conjugacy classes of newforms of level $49$ and weight $6$ that seem to be supercuspidal at $7$ the residual mod $11$ representation is never bad-dihedral. On the other hand, since $7$ is split in $F$, the restriction to $G_F$ of the projectivization of this mod $11$ representation contains a dihedral group of order $8$. Then, we can apply Lemma \ref{teo:lema3}
and conclude that the restriction to $F(\zeta_{11})$ of this residual representation is absolutely irreducible.\\
We also want to check that the mod $11$ representation of $h_{10}$ satisfies the technical condition needed to apply Theorem \ref{teo:superKisin}. We apply again Lemma \ref{teo:caruso}: since $h_2$ has weight $3$ and the $11$-adic representation is crystalline over an extension of $\Q_{11}$ of degree $2$, $2$ is even, and $(3-1)\cdot 2 < 11-1$, we see that the technical condition is satisfied. \\
Thus, if we take $h_{11}$ a minimal modular lift of this mod $11$ representation, it corresponds to a modular form of level $49$, supercuspidal at $7$, of some even weight $k \leq 12$, and we know that for the restriction to $G_F$ modularity propagates well from $h_{11}$ to $h_{10}$. To complete the proof, it suffices to show that any such $h_{11}$ can be lifted to $F$.\\
As we already mentioned, there are only two conjugacy classes of newforms in each of the spaces $S_k(49)$ with $k=4,6,8,10,12$ which are supercuspidal at $7$, and one is a twist of the other. The fields of coefficients of these newforms have degrees $1,2,4,5,6$, respectively. Since modularity is preserved by Galois conjugation and by twisting, it is enough to show that, for each $k$, one of these newforms can be lifted to $F$. We do a few computations and we observe that all these newforms have residual mod $3$ representation defined over $\F_3$, thus with image contained in the solvable group $\GL(2, \F_3)$. Moreover, in each conjugacy class there is a newform such that the $3$-adic Galois representation is ordinary, because $3 \nmid a_3$. This is the one that we will show that can be lifted to $F$. Because of the MGD prime $7$ being split in $F$ and a square modulo $3$, we also know that the mod $3$ representation, even when restricted to $G_F$, is irreducible and not bad-dihedral. Thus, the restriction to $G_F$ of the residual mod $3$ representation has irreducible solvable image and is therefore modular (as in Wiles' work [W] on modularity of elliptic curves we rely on results of Langlands and Tunnell). Since the $3$-adic representation is ordinary, the action of tame inertia on this residual representation will be given by the characters $\{\chi , 1 \}$, but it can have Serre's weight $2$ or $4$. In any case, by the M.L.T. in [SW] modularity over $F$ of this $3$-adic representation follows: we know that it is ordinary and residually modular, and in our case we can see that we can take a lift of the mod $3$ representation corresponding to a Hilbert modular form $h_{12}$ of parallel weight $2$ or $4$ and level prime to $3$, and in both cases the $3$-adic representations of $h_{12}$ will be ordinary; this is known to follow in both cases from residual ordinarity  since $3$ is split in $F$ and the $3$-adic representations are crystalline and of ``weight" $k \leq  3+1$ (in the weight $4$ case, this result is proved in [BLZ])\footnote{alternatively, see [El],[T] and [M] for similar uses of Langlands-Tunnell and Skinner-Wiles to deduce modularity of ordinary $3$-adic Galois representations of totally real number fields with solvable (irreducible) residual image}. \\
This was the last cha\^{i}non in our chain: we have seen that over $F$ modularity propagates well, starting at $h_{12}$, from any Galois representation of $G_F$ in our chain to the previous one, and this shows that the given $f$ can be lifted to $F$. 

 \section{Elementary consequences}
 
 In this section we will discuss some elementary consequences of Theorem \ref{teo:main} combined with some recent modularity results over $\Q$.
  To simplify the statements, we assume that $5$ is split in $F$.\\
 \begin{coro} \label{teo:corolario} Let $F$ be a totally real Galois number field such that the primes $2,3,5,7$ and $11$ are split in $F$. 
 Let $p$ be an odd prime. Let $$\rho_p :  G_F  \rightarrow \GL(2, \bar{\Q}_p) $$
 be a totally odd, continuous representation, with finite ramification set not containing $2$ and de Rham locally at places above $p$. Suppose that the residual representation $\bar{\rho}_p$ has non-solvable image, and that the representation $\rho_p$ can be extended to a two-dimensional Galois representation of $G_\Q$. Suppose furthermore that one of the following conditions holds:
\begin{enumerate}
	\item $\rho_p$ is, locally at places above $p$, of Hodge-Tate weights $\{0,  1   \}$ 
	\item $\rho_p$ is, locally at places above $p$, of parallel Hodge-Tate weights $\{0,  k-1  \}$, $k>2$,  and for some $v \mid p$ in $F$ the residual representation $\bar{\rho}_p$ locally at $v$ satisfies the technical condition in Theorem \ref{teo:superKisin}
\end{enumerate}
Then, the representation $\rho_p$ is modular.
 \end{coro}
 \bf{Proof}:  \rm  The proof is quite elementary. Let $\rho'_p$ be an extension of $\rho_p$ to $G_\Q$. Then, since Serre's conjecture over $\Q$ is now a Theorem (cf. [D] and [KW1]), we know that this representation is residually modular. Moreover, we can apply over $\Q$ the M.L.T. in [KiBT] and [SW] in the case of condition (1) and the one in [KiFM] in the case of condition (2) and conclude that $\rho'_p$ is modular, thus attached to a modular form $f$ whose level is odd because $\rho'_p$ is unramified at $2$, as follows from the assumptions: $\rho_p$ unramified at $2$ and $2$ split in $F$. Then we apply Theorem \ref{teo:main} and conclude that the restriction of $\rho'_p$ to $G_F$ is modular, but this restriction is precisely $\rho_p$, so this proves the corollary.
 \qed\\
 \newline
 We can also conclude modularity of a $2$-dimensional $p$-adic representation $\rho$ of $G_F$, under the assumption that $\rho$ is Galois invariant, i.e., isomorphic to all of its inner Galois conjugates. It is known that under such an assumption a suitable twist $\rho \otimes \psi$ of the representation can be extended to a $2$-dimensional representation $\rho'$ of $G_\Q$. This is proved in [Wi], section 2.4: see Lemma 5 (first assertion) for the construction of a projective representation of $G_\Q$ extending $\rho$, and Lemma 7 for the conclusion (observe that since $2$ is split in $F$ the proofs of these two Lemmas imply that $\rho'$ will have odd conductor). Finally, Lemma 6 in [Wi], under the extra assumption that $p$ is unramified in $F$, implies that if $\rho$ satisfies condition (2) in Corollary \ref{teo:corolario}, then $\rho \otimes \psi$ and $\rho'$ will also satisfy this condition. \\
 From this we can show (under the assumptions in Corollary \ref{teo:corolario}) that  $\rho \otimes \psi$ is modular, thus that $\rho$ is modular. 
 
 \begin{coro} \label{corolario2}
 Let $p$, $F$, and $\rho_p$ be as in the previous Corollary, and we keep the assumptions in the previous Corollary except that instead of assuming that $\rho_p$ can be extended to $G_\Q$ we just assume that $\rho_p$ is Galois invariant. Suppose also that $p$ is unramified in $F$.\\
 Then, $\rho_p$ is modular.
 \end{coro}
 We finish with an elementary corollary of Theorem \ref{teo:main} for the case of residual representations.
 
 \begin{coro}  \label{corolario3}
  Let $F$ be a totally real Galois number field such that the primes $2,3,5,7$ and $11$ are split in $F$. 
 Let $p$ be an odd prime. Let $$\bar{\rho}_p :  G_F  \rightarrow \GL(2, \F_{p^r}) $$
 be a totally odd representation with ramification set not containing $2$. Suppose that it is absolutely irreducible, and that it can be extended to a two-dimensional Galois representation of $G_\Q$. \\
 Then $\bar{\rho}_p$ is modular, i.e., there exists a Hilbert modular form $h$ over $F$ such that one of the $p$-adic Galois representations attached to $h$ has residual representation isomorphic to $\bar{\rho}_p$.
 \end{coro}
 If the residual image of $\bar{\rho}_p$ is solvable, modularity follows from the results of Langlands and Tunnell. If not, the proof is similar to the proof of Corollary \ref{teo:corolario} and is left to the reader.

\section{Bibliography}
[AC] Arthur, J., Clozel, L., {\it Simple algebras, base change, and the advanced theory of the trace formula}, Annals of Math. Studies {\bf 120}, Princeton U. P. (1989)
\newline
[BLGG]  Barnet-Lamb, T., Gee, T.,  Geraghty, D., {\it Congruences between Hilbert modular forms: constructing ordinary lifts}, preprint; available at www.arxiv.org 
\newline
[BLGGT]   Barnet-Lamb, T., Gee, T., Geraghty, D., Taylor, R., {\it
Potential automorphy and change of weight}, preprint; available at www.arxiv.org 
\newline
[BLZ] Berger, L., Li, H., Zhu, J., {\it Construction of some families of $2$-dimensional crystalline representations}, 
Math. Annalen {\bf 329} (2004) 365-377
\newline
[C] Caruso, X., {\it Repr\'{e}sentations semi-stables de torsion
 dans le cas $er < p - 1$}, J. Reine Angew. Math. {\bf 594} (2006) 35-92
\newline
[D] Dieulefait, L., {\it Remarks on Serre's modularity conjecture}, preprint, (2006); available at www.arxiv.org
\newline
[E] Edixhoven, B., {\it The weight in Serre's conjectures on modular forms}, Invent. Math. {\bf 109} (1992)  563-594
\newline
[El] Ellenberg, J., {\it Serre's conjecture over $\F_9$}, Annals of Math. {\bf 161} (2005)  1111-1142
\newline
[G] Geraghty, D., {\it Modularity lifting theorems for ordinary Galois
representations}, preprint
\newline
[H] Hida, H., {\it   Serre's conjecture and base change for $\GL(2)$}, Pure and Applied Math. Quarterly {\bf 5}  (2009) 81-125 
\newline
[HM] Hida, H., Maeda, Y., {\it Non-abelian Base Change for totally real fields}, Olga Taussky Todd memorial issue, Pacific J. Math. {\bf 181} (1997) 189-217 
\newline
[Kh] Khare, C., {\it A local analysis of congruences in the $(p,p)$ case: Part II}, Invent. Math. {\bf 143} (2001) 129-155 
\newline
[KW1] Khare, C., Wintenberger, J-P., {\it Serre's modularity conjecture (I)}, Invent. Math.  {\bf 178} (2009) 485-504
\newline
[KW2] Khare, C., Wintenberger, J-P., {\it Serre's modularity conjecture (II)}, Invent. Math. {\bf 178} (2009) 505-586
\newline
[KiBT] Kisin, M., {\it Moduli of finite flat group schemes, and modularity}, Annals of Math. {\bf 170} (2009) 1085-1180
\newline
[KiFM] Kisin, M., {\it The Fontaine-Mazur conjecture for $\GL_2$}, J.A.M.S {\bf 22} (2009) 641-690
\newline
[Ki2] Kisin, M., {\it  Modularity of $2$-adic Barsotti-Tate representations}, Invent. Math. {\bf 178} (2009) 587-634
\newline 
[L] Langlands, R., {\it Base Change for $\GL(2)$}, Annals of Math. Studies {\bf 96}, Princeton U. P. (1980)
\newline
[M] Manoharmayum, J., {\it On the modularity of certain $\GL(2, \F_7)$ Galois representations}, Math. Res. Letters {\bf 8} (2001) 703-712
\newline
[R1] Ribet, K., {\it On l-adic representations attached to modular forms. II}, Glasgow Math. J. {\bf 27} (1985) 185-194
\newline
[R2] Ribet, K., {\it Raising the levels of modular representations}, S\'{e}minaire de Th\'{e}orie des Nombres, Paris 1987--88, 259-271, Progr. Math. {\bf 81}, Birkhauser (1990)
\newline
[RS] Ribet, K., Stein, W., {\it Lectures on Serre's conjectures}, Arithmetic Algebraic Geometry (Park City, UT, 1999), 143-232, IAS/Park City Math. Ser. {\bf 9}, A. M. S. (2001)
\newline
[S] Serre, J-P., {\it Sur les repr{\'e}sentations modulaires de degr{\'e}
$2$ de $\Gal(\bar{\mathbb{Q}} / \mathbb{Q})$}, Duke Math. J. {\bf 54}
(1987) 179-230
\newline
[SW] Skinner, C., Wiles, A., {\it Nearly ordinary deformations of irreducible
residual representations}, Ann. Fac. Sci. Toulouse Math. (6) {\bf 10} (2001) 185-215
\newline
[T] Taylor, R., {\it On icosahedral Artin representations. II}, 
   American J. Math. {\bf 125} (2003) 549-566
   \newline
   [W] Wiles, A., {\it Modular elliptic curves and Fermat's Last Theorem}, Annals of Math. {\bf 141} (1995) 443-551
   \newline
   [Wi] Wintenberger, J.-P., {\it 
   Sur les repr{\'e}sentations $p$-adiques g{\'e}om{\'e}triques de
conducteur $1$ et de dimension $2$ de $G_\Q$}, preprint; available at www.arxiv.org 


\end{document}